\documentclass[leqno,12pt]{amsart}
\usepackage{latexsym}
\usepackage{mathrsfs}
\usepackage{amsmath}
\usepackage{amssymb}
\usepackage[bookmarks]{hyperref}
\usepackage{pstricks,pst-plot}

\font\tenscr=rsfs10 
\font\sevenscr=rsfs7 
\font\fivescr=rsfs5 
\skewchar\tenscr='177 \skewchar\sevenscr='177 \skewchar\fivescr='177
\newfam\scrfam \textfont\scrfam=\tenscr \scriptfont\scrfam=\sevenscr
\scriptscriptfont\scrfam=\fivescr
\def\scr{\fam\scrfam}

\newtheorem{theorem}{Theorem}
\newtheorem{lemma}[theorem]{Lemma}
\newtheorem{corollary}[theorem]{Corollary}
\newtheorem{proposition}[theorem]{Proposition}

\theoremstyle{definition}

\newtheorem{definition}[theorem]{Definition}

\newtheorem{example}[theorem]{Example}

\allowdisplaybreaks[1]

\newcommand{\R}{\mathbb{R}}

\newcommand{\Rd}{\R_d}
\newcommand{\CbX}{C_b(X)}

\newcommand{\bthm}{\begin{theorem}}
\newcommand{\ethm}{\end{theorem}}
\newcommand{\blem}{\begin{lemma}}
\newcommand{\elem}{\end{lemma}}
\newcommand{\bcor}{\begin{corollary}}
\newcommand{\ecor}{\end{corollary}}
\newcommand{\bprop}{\begin{proposition}}
\newcommand{\eprop}{\end{proposition}}
\newcommand{\bdefn}{\begin{definition}}
\newcommand{\edefn}{\end{definition}}
\newcommand{\bpf}{\begin{proof}}
\newcommand{\epf}{\end{proof}}

\newcommand{\ol}{\overline}

\newcommand{\bi}{\begin{itemize}}
\newcommand{\ei}{\end{itemize}}

\newcommand{\bc}{\begin{cases}}
\newcommand{\ec}{\end{cases}}

\newcommand{\ba}{\begin{array}}
\newcommand{\ea}{\end{array}}

\newcommand{\be}{\begin{equation}}
\newcommand{\ee}{\end{equation}}

\newcommand{\bea}{\begin{eqnarray}}
\newcommand{\eea}{\end{eqnarray}}

\newcommand{\beaa}{\begin{eqnarray*}}
\newcommand{\eeaa}{\end{eqnarray*}}

\newcommand{\beastar}{\begin{eqnarray*}}
\newcommand{\eeastar}{\end{eqnarray*}}

\font\tenscr=rsfs10  scaled \magstep1

\font\sevenscr=rsfs7  scaled \magstep1
\font\fivescr=rsfs5  scaled \magstep1
\skewchar\tenscr='177 \skewchar\sevenscr='177 \skewchar\fivescr='177
\newfam\scrfam \textfont\scrfam=\tenscr \scriptfont\scrfam=\sevenscr
\scriptscriptfont\scrfam=\fivescr
\def\scr{\fam\scrfam}

\def\sC{{\scr C}}

\def\vpa{\varphi_\alpha}

\def\sm{\setminus}

\begin{document}
\title[not a Borel set]{The set of bounded continuous\\ nowhere locally uniformly continuous functions is not Borel}
\author{Alexander J. Izzo}
\address{Department of Mathematics and Statistics, Bowling Green State University, Bowling Green, OH 43403}
\email{aizzo@bgsu.edu}
\thanks{The author was supported by NSF Grant DMS-1856010.}


\begin{abstract}
It is known that for $X$ a nowhere locally compact metric space, the set of bounded continuous, nowhere locally uniformly continuous real-valued functions on $X$ contains a dense $G_\delta$ set in the space $\CbX$ of all bounded continuous real-valued functions on $X$ in the supremum norm.  Furthermore, when $X$ is separable, the set of bounded continuous, nowhere locally uniformly continuous real-valued functions on $X$ is itself a $G_\delta$ set.  We show that in contrast, when $X$ is nonseparable, this set of functions is not even a Borel set.
\end{abstract}
\maketitle



\setcounter{section}{1}
We call a function $f:X\rightarrow Y$ between metric spaces \emph{locally uniformly continuous at a point $x$} if there is a neighborhood $U$ of $x$ on which $f$ is uniformly continuous.  If $f$ is locally uniformly continuous at every point of $X$, we say that $f$ is \emph{locally uniformly continuous}.  If $f$ is locally uniformly continuous at no point of $X$, we say that $f$ is \emph{nowhere locally uniformly continuous}.  Equivalently, $f$ is nowhere locally uniformly continuous if it is uniformly continuous on no open set of $X$.  

In \cite{Izzo1994} the author proved that for $X$ a separable metric space that is nowhere locally compact (i.e., locally compact at no point), the set of bounded continuous, nowhere locally uniformly continuous 
real-valued functions on $X$ is a dense $G_\delta$ set in the space $\CbX$ of bounded continuous real-valued functions on $X$ with the supremum norm.  In \cite{Izzo1999} the author extended this result by showing that in the absence of the separability hypothesis, the set of bounded continuous, nowhere locally uniformly continuous 
real-valued functions on $X$ still contains a dense $G_\delta$ set in $\CbX$, and he noted that his argument in this case did not show that the set actually \emph{is} a $G_\delta$ set.  An (anonymous) referee of \cite{Izzo1999} wrote \lq\lq $\ldots$ I believe that the collection of n.\ l.\ u.\ c.\ functions *is* a dense $G$-delta.  I encourage the author to prove this 
and re-submit the paper elsewhere.\rq\rq\ The purpose of the present paper is to show that on the contrary, for $X$ nonseparable (and nowhere locally compact) 
the set of nowhere locally uniformly continuous functions in $\CbX$ is not even a Borel set.

The proofs given here use an idea from the paper of Roberts \cite{Roberts1948} (which came to the author's attention from a footnote in the classic 
book \cite{Hurewicz-Wallman1948}).  A well-known result in dimension theory asserts that if $X$ is a separable metric space of topological dimension $n$, then in the space $I_{2n+1}^X$ of all continuous mappings of $X$ into the $(2n+1)$-dimensional Euclidean cube $I_{2n+1}$ with the uniform metric, the set of embeddings contains a dense $G_\delta$ set, and in the case when $X$ is compact, the set of embeddings \emph{is} a $G_\delta$ set.  Roberts proved that for a certain (noncompact) separable metric space $X$ the set of embeddings is \emph{not} a $G_\delta$ set in $I_{2n+1}^X$ by finding an embedded Cantor set in $I_{2n+1}^X$ whose intersection with the the set of embeddings of $X$ into $I_{2n+1}$ is not a $G_\delta$ set in the Cantor set.  We will use a slight modification of Roberts' idea.  If $Y$ and $Z$ are topological spaces and $F:Y\rightarrow Z$ is a continuous map, then $F^{-1}(E)$ is a Borel set in $Y$ for every Borel set $E$ in $Z$.  Thus to show that a set $E$ in $Z$ is not a Borel set, it is sufficient to find a continuous map $F:Y\rightarrow Z$ such that $F^{-1}(E)$ is not Borel in $Y$.

The precise statement of our result is as follows. 

\bthm\label{maintheorem}
Let $X$ be a nonseparable, nowhere locally compact metric space. 
Then the set of bounded continuous, nowhere locally uniformly continuous real-valued functions on $X$ is not a Borel set in $\CbX$.
\ethm

Before proving the theorem we present an example that illustrates the main idea behind the proof.  The reader who wishes, can skip the example and proceed directly to the proof of the theorem.  On the other hand, the reader who is content to verify only the \emph{existence} of a metric space for which the conclusion of the theorem holds can read only the example and omit the proof of the general theorem.

Throughout $\R$ will denote the real line with the standard metric.

\begin{example}
Let $S$ be a nowhere locally compact metric space with metric $d_S$, let $\Rd$ be the real line with the metric $d_d$ that takes only the values 0 and 1, and let $X=S\times\Rd$ with the metric $d$ defined by $d\bigl((s_1,\alpha_1),(s_2,\alpha_2)\bigr)=d_S(s_1,s_2) + d_d(\alpha_1,\alpha_2)$.  
We wish to show that 
the set of bounded continuous, nowhere locally uniformly continuous real-valued functions on $X$ is not a Borel set in $\CbX$.
As noted above, it will suffice to construct a continuous map $F:\R\rightarrow \CbX$ such that the inverse image under $F$ of the nowhere locally uniformly continuous functions in $\CbX$ is not Borel in $\R$. 

Fix a bounded subset $A$ of $\R$ that is not a Borel set.  Fix a bounded continuous, nowhere locally uniformly continuous real-valued function $f$ on $S$ (which exists by \cite{Izzo1999}).

For each $r\in\R$ define a function $f_r: X=S\times \Rd\rightarrow \R$ by
$$f_r(s,\alpha)= \bc
f(s) & \mbox{ if } \alpha\in\Rd\sm A\\
|\alpha-r| f(s) & \mbox{ if } \alpha\in A.\ec$$
Each function $f_r$ is bounded, and since a function on $X$ is continuous if its restriction to each set $S\times \{\alpha\}$, $\alpha\in \Rd$, is continuous, each $f_r$ is continuous.  Now define the map $F:\R\rightarrow \CbX$ by $F(r)=f_r$.  A trivial computation shows that for $p$ and $q$ in $\R$
$$\|F(p)-F(q)\|_\infty \leq |p-q| \, \|f\|_\infty,$$
where $\|\cdot\|_\infty$ denotes the supremum norm.
Thus $F$ is continuous.

Because $f$ is nowhere locally uniformly continuous, $f_r$ is nowhere locally uniformly continuous for each $r\in \R\sm A$.  In contrast, for $r\in A$ the function $f_r$ is identically zero on the open set $S\times \{r\}$, and hence, is not nowhere locally uniformly continuous.  Thus the inverse image of the bounded continuous, nowhere locally uniformly continuous functions under $F$ is the non-Borel set $\R\sm A$.  
\hfill$\square$
\end{example}

In preparation for the proof of Theorem~\ref{maintheorem}, we mention an elementary fact from descriptive set theory that will be used: There exists a bounded subset of $\R$ that is not a Borel set and whose cardinality is $\aleph_1$ (the first uncountable cardinal).  To see this, first note that if the continuum hypothesis holds, then 
every non-Borel subset of $\R$ has cardinality $\aleph_1$.  If the continuum hypothesis fails, then every subset of $\R$ of cardinality $\aleph_1$ is non-Borel, for each Borel subset of $\R$ is either countable or has the cardinality of $\R$ (see
\cite[Section 2C]{Mos} for instance).

\bpf[Proof of Theorem~\ref{maintheorem}]
We will construct a continuous map $F:\R\rightarrow \CbX$ such that the inverse image under $F$ of the nowhere locally uniformly continuous functions in $\CbX$ is not Borel in $\R$.  As noted earlier, this will establish the theorem.

Fix a bounded subset $A$ of $\R$ that is not a Borel set and has cardinality $\aleph_1$.  Because $X$ is a nonseparable metric space, there is a locally discrete collection $\sC$ of open sets of $X$ that has cardinality $\aleph_1$.  This follows from the Bing metrization theorem (and can also be proven directly).
Choose a bijection $A\rightarrow \sC$ and denote the member of $\sC$ corresponding to $\alpha\in A$ by $U_\alpha$.  For each $\alpha$, choose a nonempty open set $V_\alpha$ such that $\ol V_\alpha\subset U_\alpha$.  Then choose, for each $\alpha$, a locally uniformly continuous function $\vpa:X\rightarrow [0,1]$ such that $\vpa(\ol V_\alpha)=\{1\}$ and $\vpa(X\sm U_\alpha)=\{0\}$.  For instance, one can define $\vpa$ by
\begin{equation}\label{vpa}
\vpa(x)=\frac{d_{X\sm U_\alpha}(x)}{d_{X\sm U_\alpha}(x) +d_{V_\alpha}(x)}
\end{equation}
where $d_E(x)=\inf\{d(x,e):e\in E\}$ (the distance from the point $x$ to the set $E$).
The easy proof that the function $\vpa$ defined by equation~(\ref{vpa}) is locally uniformly continuous is given in \cite[Lemma~5]{Izzo1994}.  Choose a bounded continuous, nowhere locally uniformly continuous real-valued function $f$ on $X$ (which exists by \cite{Izzo1999}).

Now for each $r\in\R$, define $f_r:X\rightarrow \R$ by
$$f_r(x)=\biggl[ 1-\sum_{\alpha\in A} \bigl(1-|\alpha - r | \bigr) \vpa(x) \biggr ] f(x).$$
Because each point of $X$ has a neighborhood on which at most one of the $\vpa$ is not identically zero, $f_r$ is a well-defined continuous function.  Note also that $f_r$ is a bounded function.  Now define the map $F:\R\rightarrow \CbX$ by $F(r)=f_r$.

For $p$ and $q$ in $\R$ and $x$ in $X$ we have

\begin{equation*}
\begin{split}
\bigl | f_p(x) - f_q(x) \bigr |& = \Biggl | \biggl [\, \sum_{\alpha\in A} \Bigl ( |\alpha - q | - |\alpha - p | \Bigr ) \vpa(x) \biggr ] f(x) \Biggr |\\ 
& \leq \Biggl |  \sum_{\alpha\in A} \Bigl ( |\alpha - q | - |\alpha - p | \Bigr ) \vpa(x) \Biggr | \|f\|_\infty\\ 
& \leq \biggl [\, \sum_{\alpha\in A} | p-q | \vpa(x) \biggr ] \|f\|_\infty\\ 
&\leq | p-q | \, \|f\|_\infty
\end{split}
\end{equation*}
since $\vpa(x)$ is nonzero for at most one $\alpha$.  Thus $F$ is continuous.

To complete the proof we show that the inverse image of the set of bounded continuous, nowhere locally uniformly continuous functions under $F$ is the non-Borel set $\R\sm A$.  Observe that for $r\in A$, the function $f_r$ is identically zero on $V_r$, and hence, is not nowhere locally uniformly continuous.  Now consider $r\in \R\sm A$.  Define $g_r:X\rightarrow \R$ by
$$g_r(x)= 1-\sum_{\alpha\in A} \bigl(1-|\alpha - r | \bigr) \vpa(x)$$
so that $f_r=g_rf$.  For each $x\in X$, each term in the sum over $A$ in the definition of $g_r$ is strictly less than 1, and since at most one of the terms is nonzero, $g_r$ is zero free.  Furthermore, since each $\vpa$ is locally uniformly continuous and the supports of the $\vpa$ form a locally discrete family, $g_r$ is locally uniformly continuous.  Thus the following lemma yields that $f_r$ is nowhere locally uniformly continuous thereby concluding the proof.
\epf

\blem
The product of a zero free locally uniformly continuous real-valued function and a continuous, nowhere locally uniformly continuous real-valued function is nowhere locally uniformly continuous.
\elem

\bpf
Suppose that $f$, $g$, and, $h$ are real-valued functions such that $h=gf$, that $g$ is zero free and locally uniformly continuous, and that $f$ is continuous.  We show that if $h$ fails to be nowhere locally uniformly continuous, then so does $f$.

Let $U$ be an open set on which $h$ is uniformly continuous.  By shrinking $U$, we can assume that $h$ is bounded on $U$.  By shrinking $U$ further, we can also assume that $g$ is uniformly continuous on $U$ and that $g$ is bounded away from zero on $U$.  Then $1/g$ is bounded and uniformly continuous on $U$.  Thus the restriction of $f=h/g$ to $U$ is a product of two bounded uniformly continuous functions and hence is uniformly continuous.
\epf

\end{document}